# Computing persistent homology by spanning trees and critical simplices

Dinghua Shi[1], Zhifeng Chen[2], Chuang Ma[3], Guanrong Chen[4]

[1]Department of Mathematics, College of Science, Shanghai University, Shanghai, China.
[2]School of Big Data, Fuzhou University of International Studies and Trade, Fuzhou, China.
[3] Department of Data Science and Big Data Technology, School of Internet, Anhui University, Hefei, China.
[4]Department of Electrical Engineering, City University of Hong Kong, Hong Kong, China.
Emails: shidh2012@sina.com; 920978196@qq.com; chuang_m@126.com; eegchen@cityu.edu.hk

**Abstract:** Topological data analysis can extract effective information from higher-dimensional data. Its mathematical basis is persistent homology. The persistent homology can calculate topological features at different spatiotemporal scales of the dataset; that is, establishing the integrated taxonomic relation among points, lines and simplices. Here, the simplicial network composed of all-order simplices in a simplicial complex is essential. Because the sequence of nested simplicial subnetworks can be regarded as a discrete Morse function from the simplicial network to real values, a method based on the concept of critical simplices can be developed by searching all-order spanning trees. Employing this new method, not only the Morse function values with the theoretical minimum number of critical simplices can be obtained, but also the Betti numbers and composition of all-order cavities in the simplicial network can be calculated quickly. Finally, this method is used to analyze some examples and compared with other methods, showing its effectiveness and feasibility.

## Introduction

In the era of bigdata[1], higher-dimensional data and their analysis connect together the fields of data science[2-3], network science[4] and computational topology[5-6]. In these studies, algebraic topology especially persistent homology[7-10] plays a key role. In fact, the shape of a dataset can be viewed as the real representation of the true data only if it appears persistently in various spatiotemporal scales; otherwise, it might likely contain sampling errors or noise.

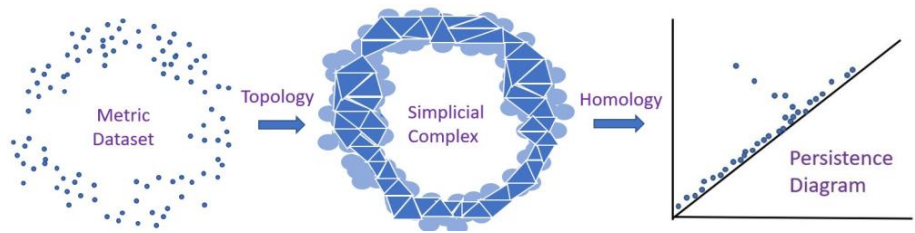

Fig. 1. Point-cloud data and their topological features (Ref. [9])

(left, point-cloud data; middle, simplicial complex; right, persistence diagram)

One typical example is the point-cloud data, such as the recorded spatial coordinates and brightness of a scanned image, as shown in Fig. 1 (left picture). With this dataset, with a given threshold value, one can establish the topological relationship among all the points, obtaining a simplicial complex as shown in Fig. 1 (middle picture). Simple complexes obtained by different methods are usually different; for example, they can be dense such as Čech complex[11] and Vietoris-Rips complex[12], or sparse such as Alpha complex[13] and Witness complex[14]. From the perspective of network science, the network composed of all-order (all different orders of) simplices in a simplicial complex is called a simplicial network. To this end, by increasing the threshold values, one can obtain and analyze the filtered simplicial network so as to dynamically extract its topological features under different scales, as shown in Fig. 1 (right picture). In general, time series data can be converted to point-cloud data after determining the number of points, time windows and spatial dimensions[15].

Network science can be traced back to Euler when he solved the Königsberg seven-bridge problem, thereby established the graph theory. In the era of the Internet with bigdata, the models of small-world networks[16] and scale-free networks[17] attracted a lot of attention. Consequently, in the study of totally homogeneous networks[18], it was found that the Euler characteristic number $\chi$ can be applied to analyzing complex networks, where cliques (fully-connected subnetworks) are simplices of different orders, such as node (0-simplex, i.e., order 0), edge (1-simplex), triangle (2-simplex), and so on. The summation of these clique numbers with alternative signs is an invariant value of the given network, namely the Euler characteristic number $\chi$. On the other hand, in addition to cliques, there are many cavities in a large-scale network. The shape and the number of some cavities commonly seen in geometry are shown by the examples in Fig. 2. In order to distinguish the numbers of cavities on a sphere and on a torus, Poincaré introduced the concept of triangulation, to triangulate a cycle to get a triangle that is not a 2-simplex, to triangulate a sphere to get a tetrahedron that is not a 3-simplex, and to triangulate a torus to get a more complex network, which will be further discussed in the **Method** section. The numbers of cavities are called Betti numbers, denoted $\beta_k$, and the summation of the Betti numbers with alternative signs is equal to $\chi$ (see Fig. 2).

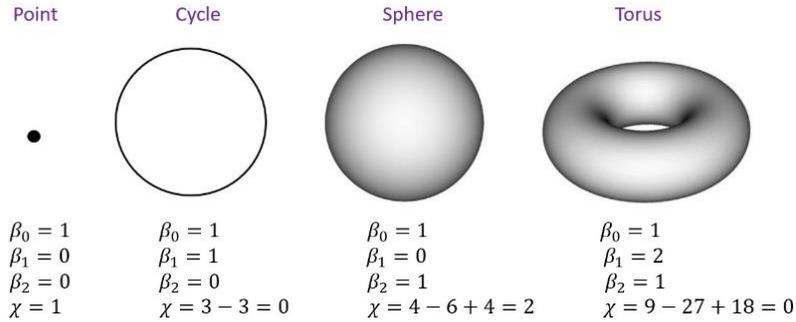

Fig. 2. Betti numbers of a node, a cycle, a sphere and a torus (Ref. [19]) and their characteristic number

One conventional approach to network science studies is to consider simplices as basic network elements and define their vector spaces, in which furthermore introduce chains and cycles etc., and consequently apply the theory and methods of algebraic topology. Here, these networks are called simplicial (complexes[19]) networks[20], which are of higher-order networks[21-22] and their orders are determined by the highest-order simplex. Thus, the common graph networks are 1-order ones because they only consider nodes and edges. It is known that the cyclic structures in simplicial networks are of fundamental importance because they provide feedback paths in higher-order dynamical interactions over such networks[23]. However, for relatively dense networks, the numbers of cycles are extremely large. Therefore, classifying cycles in large simplicial networks is a great challenge, which eventually relies on the computation of homology groups[24]. In addition to Betti numbers, which can only give the ranks of homology groups, the cycles whose classes constitute the elements of the homology groups carry some important information. In practice, interests are usually in representative cycles that have some optimal properties.

In homology theory, persistent homology starts from the Morse theory[25], which uses a continuous real variable function to calculate homology. Barannikov[26] studied framed Morse complex and its invariants. Forman[27] discredited Morse functions, and assigned real values to simplices in simplicial networks thereby classifying them. It was revealed[26-27] that the number of critical simplices (see "Basic Concepts" below for a precise definition) is between the number of cliques and the number of cavities, and that the summation of the numbers of critical simplices with alternative signs is also equal to the characteristic number $\chi$. There are several different ways to assign values to simplices; for example, Sizemore[28] used the edge weights of a weighted network, Horak[29] assigned values to simplices from lower-order to higher-order, and Kannan[30] designed a method that can theoretically find a near minimum number of critical simplices.

Similar to the point-cloud data scenario where different thresholds can induce a filtered simplicial network, by assigning values of a discrete Morse function to simplices one can also obtain a sequence of nested simplicial subnetworks[31,32]. However, different methods for assigning values will result in different filtered networks with different numbers of critical simplices. The results of computing persistent homology can be represented by persistence diagrams, persistence barcodes, persistence landscapes, and so on. However, the persistence barcodes of different Morse values are different too, as can be seen from Fig. 3 in the **Method** section.

This paper proposes an optimal method for assigning Morse values by using all-order spanning trees of the network, such that the number of critical simplices in each order is equal to the Betti number of the same order. Furthermore, by solving the systems of spanning trees and critical simplex equations, those simplices that constitute all-order cavities can be easily obtained because the corresponding equations have unique solutions.

The rest of the paper is organized as follows. In the **Method** section, homology of simplicial network and discrete Morse function will be introduced, where an example is given to show how to use spanning trees to assign Morse function values to simplicial networks, and to draw persistence barcodes. Then, computational methods based on the systems of spanning trees and critical simplex equations will be developed. In the **Results** section, the new method will be applied to some networks, specifically the C. elegans neural network, BA scale-free model network, and a network formed by point-cloud dataset, which will be compared with the existing methods to demonstrate the effectiveness and feasibility of the new method. Finally, in the **Discussion** section, different sequences of nested networks, the optimal property of the obtained cavities, and the influence of various chains on homology calculation, will be discussed.

## Method

Although homology calculation involves classifying cycles of a simplicial network, one can change the perspective to consider spanning trees of the network instead. This is because each new simplex added to a spanning tree will create a unique cycle. To ensure that the generated cycle is a cavity rather than just another simplex of higher order, the newly added simplex must be a critical one. For this purpose, an optimal method of assigning Morse values not only has to keep the persistence of cavities but also needs to reduce their calculations. It will be shown that the new method for solving the system of spanning trees and critical simplex equations can greatly improve the computational efficiency, e.g., comparing to the 0-1 programming method used for searching all cycles[24].

### Basic concepts

A $k$-simplex $\alpha_k$ is a fully-connected subnetwork composed of $k+1$ nodes, denoted by $(v_0, v_1, \cdots, v_k)$. Let $p$ and $q$ be two integers. A simplex $\alpha_p$ is a face of another simplex $\alpha_q$ if $\alpha_p \subset \alpha_q$. A simplex $\alpha_q$ is a coface of another simplex $\alpha_p$ if $\alpha_q \supset \alpha_p$. A simplicial complex $K$ satisfies that (i) its every node is in some simplex in $K$, and (ii) for a simplex $\alpha \in K$, if a simplex $\beta \subset \alpha$ then $\beta \in K$. A conventional network is denoted by $G = \{V, E\}$, where $V$ is the node set and $E$ is the edge set, which has an adjacency matrix to characterize the connected edges. When all-order simplices of a network are enumerated or calculated, the network is called a simplicial network and denoted as $K = \{V, E, T, \cdots\}$, which has many incidence matrixes, for example boundary matrixes $B_1$, $B_2$, etc., to describe the relationship between simplices[20]. A filtration of a simplicial network is defined as a sequence of nested simplicial subnetworks, $\varnothing \subseteq K_0 \subseteq K_1 \subseteq \cdots \subseteq K_n = K$.

Let $C_k$ be the vector space in the binary field with a basis consisting of $k$-simplices, whose dimension $m_k$ is equal to the number of $k$-simplices. In the binary field, the addition between two vectors $c$ and $d$ is defined by set operations as $c+d = (c \cup d)-(c \cap d)$. To study the linear dependence, boundary matrices are introduced. For instance, in vector space $C_1$, define a node-edge matrix $B_1$, in which an entry is 1 if a node is in the corresponding edge; otherwise, it is 0. The rank $r_k$ ($r_0 = 0$ by convention) of the boundary matrix $B_k$ is the number of linearly independent vectors in space $C_k$. Furthermore, define a boundary operator $\partial_k$: $C_k \to C_{k-1}$ to connect two successive spaces. Linear combinations of elements in $C_k$ are called $k$-chains. A $k$-cycle $l$ is defined by $\partial_k(l) = 0$. Moreover, define the kernel space of $C_k$ by $\ker(\partial_k) = \{l \in C_k \mid \partial_k(l) = 0\}$, and denote it as $Z_k$. Also, define the image space on $C_k$ by $\text{im}(\partial_{k+1}) = \{\partial_{k+1}(l) \mid l \in C_{k+1}\}$, which is the image of $C_{k+1}$ mapping to $C_k$ and denote it as $Y_k$. Since $\partial_k(\partial_{k+1}) = 0$, one has $\text{im}(\partial_{k+1}) \subseteq \ker(\partial_k)$. Two $k$-cycles $c$ and $d$ are said to be equivalent, denoted $c \sim d$, if $c+d$ is a boundary of a $(k+1)$-chain. Classifying the kernel (cycle) space $Z_k$ with respect to $Y_k$ yields a homology group $Z_k/Y_k$. A $k$-cavity is a cycle that is usually selected with the shortest length in one of the linearly independent cycle-equivalent classes, whose number is equal to the Betti number $\beta_k = m_k - r_k - r_{k+1}$.

A simplicial network containing no $k$-cycle is called a forest, while a connected forest is called a $k$-tree. If the network itself is not a $k$-tree, but it can be seen as a certain $k$-tree with some additional simplices of older not higher than $k$, then this $k$-tree is called a $k$-order spanning tree. In this paper, the $k$-order spanning tree is searched by finding the rank of the boundary matrix $B_k$, for which details given in **Supplementary Table S1**.

Now, let $K$ be a simplicial network. Given a function $f$: $K \to R$, for each simplex $\alpha_p \in K$, define two sets of simplices: $U^f(\alpha_p) = \{\alpha_{p+1} \mid \alpha_{p+1} \supset \alpha_p, f(\alpha_{p+1}) \leq f(\alpha_p)\}$ and $V^f(\alpha_p) = \{\mid \alpha_{p-1} \subset \alpha_p, f(\alpha_{p-1}) \geq f(\alpha_p)\}$. Then, a function $f$: $K \to R$ is a discrete Morse function if and only if for every simplex $\alpha_p \in K$, both $\#U^f(\alpha_p) \leq 1$ and $\#V^f(\alpha_p) \leq 1$ hold. The symbol # indicates the number of elements in the set. That is, at most one coface of one order higher is allowed such that $f(\alpha_{p+1}) \leq f(\alpha_p)$ and at most one face of one order lower is allowed such that $f(\alpha_{p-1}) \geq f(\alpha_p)$. A simplex $\alpha_p$ of a simplicial network $K$ with a discrete Morse function $f$ is **critical** if and only if both $\#U^f(\alpha_p) = 0$ and $\#V^f(\alpha_p) = 0$ hold. Let $c_k$ be the number of $k$-order critical simplices. Then, for an $l$-order simplicial network, one has $m_k \geq c_k \geq \beta_k$ when $0 \leq k \leq l$ and $\chi = m_0 - m_1 + m_2 - m_3 + \cdots + (-1)^l m_l = c_0 - c_1 + c_2 - c_3 + \cdots + (-1)^l c_l = \beta_0 - \beta_1 + \beta_2 - \beta_3 + \cdots + (-1)^l \beta_l$.

If the function defined on a simplicial network $K$ can induce a sequence of nested simplicial subnetworks, $\varnothing \subseteq K_0 \subseteq K_1 \subseteq \cdots \subseteq K_n = K$, then this sequence is called a filtered simplicial network. A $k$-cavity that appears in the subnetwork $K_i$ can potentially become the boundary of a $(k+1)$-chain of a later subnetwork $K_{i+j}$ with $j > 0$, which will no longer constitute a $k$-cavity in $K_{i+j}$. As such $k$-cavity has a unique index that corresponds to its birth ($i$ of $K_i$) and death ($i+j$ of $K_{i+j}$) in crossing the filtration, the other $k$-cavities without dead indexes are said to be persistent.

Now, consider the network with seven 0-simplices (nodes), ten 1-simplices (edges), and three

2-simplices (shaded triangles), as shown in Fig. 3 (a). The characteristic number of this 2-order simplicial network **K** is $\chi=7-10+3=0$. Next, specify the Morse values to the simplex for the network in the following way[32]: 1-$v_0$, 2-$v_1$, (1,2)-$e_1$, 3-$v_2$, 4-$v_3$, (3,4)-$e_3$, 5-$v_4$, (4,5)-$e_4$, 6-$v_5$, (5,6)-$e_5$, (1,2)-$e_6$, (4,6)-$e_7$, (4,5,6)-$f_7$, 7-$v_8$, (3,7)-$e_8$, (4,7)-$e_9$, (3,6)-$e_{10}$, (2,6)-$e_{11}$, (2,3,6)-$f_{11}$, (3,4,7)-$f_{12}$. Thus, the network has 2 nodes $v_0$, $v_2$; 3 edges $e_6$, $e_9$, $e_{10}$; and 1 face $f_{12}$, which are critical simplices. The summation of the numbers of critical simplices with alternative signs, $2-3+1=0$, is also equal to the characteristic number $\chi$. The sequence of nested simplicial subnetworks of this network **K** is $\emptyset \subseteq K_0 \subseteq K_1 \subseteq \cdots \subseteq K_{12} = K$. For example, $K_2 = \{v_0, v_1, e_1, v_2\}$. Non-diamond arrow lines in Fig. 3 (b) show the corresponding Betti persistence barcodes, where $v_0$ and $e_{10}$ are persistent, while $v_2$ and $e_9$ are birth and death elements. That is, $v_0$ and $e_{10}$ are cavities of order 0 and 1, respectively. The starting time is indicated by the subscript. Here, $v_2$ shows that the birth time is 2, $e_6$ shows that the death time is 6; $e_9$ shows that the birth time is 9, $f_{12}$ shows that the death time is 12.

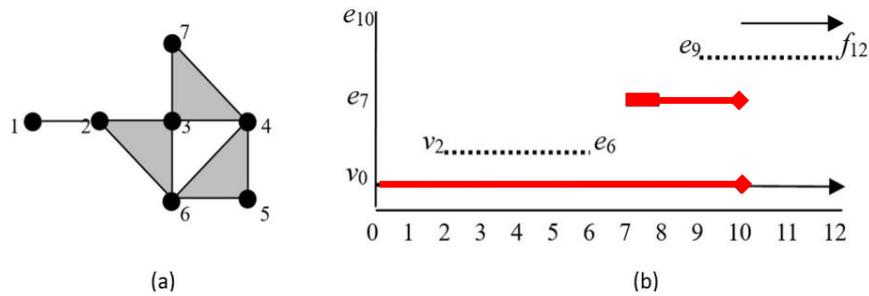

Fig. 3  (a) A 2-order simplicial network with three shaded triangles;
(b) Betti persistence barcodes of the simplicial network with different Morse values

However, different methods of assigning values may result in different starting time and also different sizes of a cavity. The **red diamond arrow** lines in Fig. 3 (b) show the Betti persistence barcodes of the corresponding Morse values: 1-$v_0$, 2-$v_1$, (1,2)-$e_1$, 3-$v_2$, (2,3)-$e_2$, 4-$v_3$, (3,4)-$e_3$, 5-$v_4$, (4,5)-$e_4$, 6-$v_5$, (5,6)-$e_5$, 7-$v_6$, (3,7)-$e_6$, (3,6)-$e_7$, (4,6)-$e_8$, (4,5,6)-$f_8$, (4,7)-$e_9$, (3,4,7)-$f_9$, (2,6)-$e_{10}$, (2,3,6)-$f_{10}$. Clearly, there are only 2 critical simplices: node $v_0$ and edge $e_7$. So the figure between critical simplices and their subscript numbers is just a persistence barcodes. The summation of the numbers of critical simplices with alternative signs, $1-1=0$, is also equal to the characteristic number $\chi$. The sequence of nested simplicial subnetworks is $\emptyset \subseteq K_0 \subseteq K_1 \subseteq \cdots \subseteq K_{10} = K$.

The above new Morse values show that the original $e_{10}$ comes before $e_7$, the 1-order cavity starting with 4 edges, but after joining by $e_8$ it has only 3 edges. This means that the 1-order cavity appears earlier and then shrinks from big to small.

**An example**

A torus with small red cycles and big blue cycles is shown in Fig. 4 (a) for illustration. Consider a network of three small red cycles and three big blue cycles, which has 9 nodes and 18 edges.

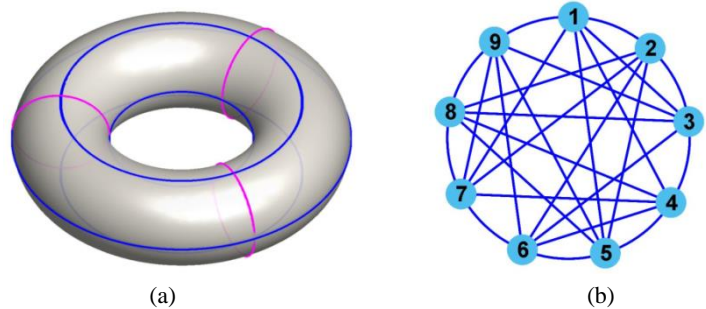

Fig. 4  (a) A torus example; (b) the obtained network with 9 nodes and 27 edges

A triangulation of the torus can be obtained by adding 9 edges to the quadrilateral in the network. The obtained network with 9 nodes and 27 edges is shown in Fig. 4 (b). The simplicial network composing of all-order simplices in the network is obtained as $K$ = {9  0-simplices (nodes), 27 1-simplices (edges), 18  2-simplices (triangles)}. Its characteristic number $\chi = m_0 - m_1 + m_2 = 9 - 27 + 18 = 0$. The boundary matrixes $B_1$, $B_2$, etc. describe the relationship between simplices. The 18 2-order simplices of $K$ are: (1,2,5), (1,2,7), (1,3,4), (1,3,9), (1,4,5), (1,7,9), (2,3,6), (2,3,8), (2,5,6), (2,7,8), (3,4,6), (3,8,9), (4,5,8), (4,6,7), (4,7,8), (5,6,9), (5,8,9), (6,7,9). However, triangles (1,2,3), (1,4,7), (1,5,9), (2,5,8), (2,6,7), (3,4,8), (3,6,9), (4,5,6), (7,8,9) are not 2-simplices in $K$.

**(1)** Searching all-order spanning trees in a network

There are many ways to find spanning trees in a given network. Here, the row elementary transformation method of matrix is used to find the rank of the boundary matrixes $B_1$ and $B_2$. Start from the last row of the boundary matrix $B_1$ to find the leftmost 1, and go up row by row. If two different rows in the same column are both 1, by the row transformation over binary operations they will be removed. The desired spanning tree is finally obtained, as can be seen from the columns corresponding to all red-colored 1 (may be an identity matrix) shown in Table 1.

The 1-order spanning tree of the triangulation network consists of (1,2), (1,3), (1,4), (1,5), (1,7), (1,9), (2,6), (2,8). This shows that $r_1 = 8$.

The 2-order spanning tree of the triangulation network consists of (1,2,5), (1,2,7), (1,3,4), (1,3,9), (1,4,5), (1,7,9), (2,3,6), (2,3,8), (2,5,6), (2,7,8), (3,4,6), (3,8,9), (4,5,8), (4,6,7), (4,7,8), (5,6,9), (5,8,9). This shows that $r_2 = 17$. Thus $\beta_0 = 9 - 0 - 8 = 1$, $\beta_1 = 27 - 8 - 17 = 2$ and $\beta_2 = 18 - 17 - 0 = 1$.

**(2)** Classifying the remaining simplices other than spanning trees

Except for the simplices in the 1-order spanning tree, the remaining 1-simplices can be classified into two categories: simplices in the 2-order spanning tree and 1-order cavity-generating simplices. For example, 1-simplices (2,3) and (4,7) are not in the 2-order spanning tree. Thus, they are 1-order cavity-generating simplices. Similarly, the 2-simplex (6,7,9) is a 2-order cavity-generating simplex. Obviously, the first selected node 1 is a 0-order cavity-generating simplex.

**(3)** Assigning Morse function values to simplices

To find the lower bound of the critical simplex numbers, it should be noted that only those cavity-generating simplices can be critical simplices. Therefore, the function values are specified to the simplices in the following way: 1-$v_0$, 2-$v_1$, (1,2)-$e_1$, 3-$v_2$, (1,3)-$e_2$, 4-$v_3$, (1,4)-$e_3$, 5-$v_4$, (1,5)-$e_4$, 7-$v_5$, (1,7)-$e_5$, 9-$v_6$, (1,9)-$e_6$, 6-$v_7$, (2,6)-$e_7$, 8-$v_8$, (2,8)-$e_8$; (2,3)-$e_9$, (4,7)-$e_{10}$; (2,5)-$e_{11}$, (1,2,5)-$f_{11}$, (2,7)-$e_{12}$, (1,2,7)-$f_{12}$, (3,4)-$e_{13}$, (1,3,4)-$f_{13}$, (3,9)-$e_{14}$, (1,3,9)-$f_{14}$, (4,5)-$e_{15}$, (1,4,5)-$f_{15}$, (7,9)-$e_{16}$, (1,7,9)-$f_{16}$, (3,6)-$e_{17}$, (2,3,6)-$f_{17}$, (3,8)-$e_{18}$, (2,3,8)-$f_{18}$, (5,6)-$e_{19}$, (2,5,6)-$f_{19}$, (7,8)-$e_{20}$, (2,7,8)-$f_{20}$; (4,6)-$e_{21}$, (3,4,6)-$f_{21}$, (8,9)-$e_{22}$, (3,8,9)-$f_{22}$, (5,8)-$e_{23}$, (4,5,8)-$f_{23}$, (6,7)-$e_{24}$, (4,6,7)-$f_{24}$, (4,8)-$e_{25}$, (4,7,8)-$f_{25}$, (6,9)-$e_{26}$, (5,6,9)-$f_{26}$, (5,9)-$e_{27}$, (5,8,9)-$f_{27}$; (6,7,9)-$f_{28}$.

The subscripts of critical simplices, 1-$v_0$, (2,3)-$e_9$, (4,7)-$e_{10}$, (6,7,9)-$f_{28}$, indicate the starting time of one 0-order cavity, two 1-order cavities and one 2-order cavity, respectively.

This Morse function induces a sequence of nested simplicial subnetworks, $\varnothing \subseteq K_0 \subseteq K_1 \subseteq \cdots \subseteq K_{28} = K$, for the triangulation network. It will greatly simplify the homology calculation for lower orders.

**(4)** Computing simplices composed of all-order cavities

A $k$-cavity can be expressed as $x = (x_1, x_2, \cdots, x_{m_k}) \in C_k$, in which each component $x_i$ takes value 1 or 0, where 1 represents a $k$-simplex with index $i$ in the cavity while 0 means no such simplex exists. Let $B_k$ be the boundary matrix between the ($k-1$)-simplices and the $k$-simplices. Then, a $k$-cavity $x$ must satisfy the equation $B_k x^T = 0 \pmod 2$.

The $j$th-column of matrix $B_k$ is marked as $B_k^j$, and the columns in the matrix obtained by removing the $j$th-column from $B_k$ is denoted as $B_k^{-j}$. Suppose that the $j$th-simplex is a $k$-order cavity-generating simplex. From the equation $B_k x^T = 0$ (mod 2), with multiplication according to block matrices, it follows that $B_k^{-j}(x^{-j})^T = -B_k^{j}$ (mod 2). The matrix composed of all columns of the $k$-order spanning tree and the matrix composed of all columns of the corresponding $k$-order cavity-generating simplices in matrix $B_k$ are denoted as $T$-$B_k$ and $C$-$B_k$, respectively. Let $T$-$x$ express the unknown 0-1 matrix of simplices corresponding to the $k$-order spanning tree. Then, all simplices of the $k$-order cavities satisfy the following matrix equation:

$$(T\text{-}B_k)(T\text{-}x)^T = (C\text{-}B_k) \text{ (mod 2)}. \tag{1}$$

Note, however, that although the spanning tree matrix ($T$-$B_k$) is of full rank, the inverse of the matrix $(T\text{-}B_k)^T(T\text{-}B_k)$ may not exist because of the relationship between matrix operations over the binary field, thus one can only use the row operation of the matrix elementary transformation to solve this matrix equation.

**Table 1** Solving matrix equation (1) by row operations over the binary field

|   | $T$-$B_1$ | | | | | | | | $(T\text{-}x)^T$ | $C$-$B_1$ | |
|---|---|---|---|---|---|---|---|---|---|---|---|
|   | (1,2) | (1,3) | (1,4) | (1,5) | (1,7) | (1,9) | (2,6) | (2,8) |  | (2,3) | (4,7) |
| 1 | 1 | 1 | 1 | 1 | 1 | 1 | 0 | 0 |  | 0 | 0 |
| 2 | 1 | 0 | 0 | 0 | 0 | 0 | 1 | 1 | $x_{(1,2)}$ | 1 | 0 |
| 3 | 0 | 1 | 0 | 0 | 0 | 0 | 0 | 0 | $x_{(1,3)}$ | 1 | 0 |
| 4 | 0 | 0 | 1 | 0 | 0 | 0 | 0 | 0 | $x_{(1,4)}$ | 0 | 1 |
| 5 | 0 | 0 | 0 | 1 | 0 | 0 | 0 | 0 | $x_{(1,5)}$ | 0 | 0 |
| 6 | 0 | 0 | 0 | 0 | 0 | 0 | 1 | 0 | $x_{(1,7)}$ | 0 | 1 |
| 7 | 0 | 0 | 0 | 0 | 1 | 0 | 0 | 0 | $x_{(1,9)}$ | 0 | 0 |
| 8 | 0 | 0 | 0 | 0 | 0 | 0 | 0 | 1 | $x_{(2,6)}$ | 0 | 0 |
| 9 | 0 | 0 | 0 | 0 | 0 | 1 | 0 | 0 | $x_{(2,8)}$ | 0 | 0 |

From Table 1, one obtains the solution of matrix equation (1): $x_{(1,2)}=(1,0)$, $x_{(1,3)}=(1,0)$, $x_{(1,4)}=(0,1)$, $x_{(1,5)}=(0,0)$, $x_{(1,7)}=(0,1)$, $x_{(1,9)}=(0,0)$, $x_{(2,6)}=(0,0)$, $x_{(2,8)}=(0,0)$. Thus, two 1-order cavities are obtained: {(1,2), (1,3), (2,3)} and {(1,4), (1,7), (4,7)}, respectively. Similarly, a 2-order cavity consists of {(1,2,5), (1,2,7), (1,3,4), (1,3,9), (1,4,5), (1,7,9), (2,3,6), (2,3,8), (2,5,6), (2,7,8), (3,4,6), (3,8,9), (4,5,8), (4,6,7), (4,7,8), (5,6,9), (5,8,9), (6,7,9)}. This 2-order cavity is obtained by similarly solving a system of equations, which also yields a triangulation of the torus.

All computations and tables of this example are presented in **Supplementary Table S1**.

**Method summary**

The method developed in this paper is an optimal method of assigning Morse values to a given simplicial network. The method **starts** from any node in a connected branch and **assigns** it a value of 0; **searches** all-order spanning trees and **identifies** critical simplices with the boundary matrix $B_k$; **finds** the node and edge connected to the start node in 1-oeder spanning tree and **assigns** a value of PO to the node and edge until all edges in the spanning tree are exhausted, then **assigns** the value of 1-oeder critical simplex one by one until the last $l$-order simplex is reached for an $l$-order simplicial network; **gets** the nested simplicial subnetworks $\varnothing \subseteq K_0 \subseteq K_1 \subseteq \cdots \subseteq K_n = K$, where $n = r_1+\beta_1+r_2+\beta_2+\cdots+\cdots+r_l+\beta_l$. Then, the method **solves** the equation systems of spanning trees and critical simplices from lower order to higher order, and finally **obtains** simplices of all $k$-order cavities by row operations on a matrix elementary transformation in the binary field.

**Results**

The procedure shown in the **Method** section is summarized as follows. Firstly, a simplicial network is obtained by calculating or enumerating all simplices. Secondly, in the simplicial network, a sequence of nested simplicial subnetworks (i.e., filtered simplicial network) is obtained by searching all-order spanning trees and assigning Morse function values to simplices. Finally, in different subnetworks, simplices composed of all-order cavities are obtained by solving matrix equation $(T\text{-}B_k)(T\text{-}x)^T = (C\text{-}B_k)$ (mod 2).

Now, the developed method is applied to the C. elegans neural network, BA scale-free network, and Stanford Dragon graphic network. The results will be compared with the existing methods.

**C. elegans neural network**[33]

As a conventional network $G = \{V, E\}$, it has 297 neurons in $V$ and 2148 synapses in $E$, and the adjacency matrix is given by the dataset from Ref. [33].

By computing the conventional network $G^{24}$, the following simplicial network $K = \{V, E, T, \cdots\}$ can be obtained, with the Euler characteristic number (and the Betti numbers $\beta_k$) $\chi = 297-2148+3241-2010+801-240+40-2 \ (= 1-139+121-4) = -21$. The ranks of boundary matrixes $B_1, B_2, \cdots, B_7$ for describing the relationship between simplices are $r_1=296, r_2=1713, r_3=1407, r_4=599, r_5=202, r_6=38, r_7=2$, respectively. Then, one can search for four 3-order cavities by 0-1 programming.

The main results (with all details included in **Supplementary Table S2**) are as follows:

(1) All-order spanning trees in the network are obtained, where the number of simplices in the $k$-order spanning tree is $r_k$, $k = 1, 2, \cdots, 7$;

(2) Morse function values and critical simplices in the network are assigned, where the numbers of the critical simplices in the $k$-simplices are $\beta_k$, $k = 0, 1, 2, 3$, and the rest are zero;

(3) Nested simplicial subnetworks $\varnothing \subseteq K_0 \subseteq K_1 \subseteq \cdots \subseteq K_{4521} = K$ in the network are obtained, where $n = r_1+\beta_1+r_2+\beta_2+r_3+\beta_3+r_4+\cdots+r_7=4521$;

(4) Simplices composed of all-order cavities in the network and length distributions of cavities are obtained.

For example, the length distribution of 1-order cavities is shown in the left side of Table 2, while the optimal length distribution of 1-order cavities is on the right side.

**Table 2** Length distribution of 1-order cavities

| length | 4  | 5  | 6  | 7  | 8  | 9 | 10 | 4   | 5 |
|--------|----|----|----|----|----|---|----|-----|---|
| number | 16 | 48 | 29 | 23 | 18 | 4 | 1  | 138 | 1 |

The optimal length distribution of 1-order cavities is obtained by the exhaustive searching method. First, according to the critical simplices, search all 1-order cycles with length of 4 that pass through each given critical simplex. Then, select those cycles that are linearly independent as cavities and delete the corresponding critical simplices. Finally, increase the lengths of the searched cycles in the remaining critical simplices and search again, until all cavities are found.

The length distribution of 2-order cavities appears to be more complicated. Tracing the evolution of cavities from large to small will find the length distribution of 2-order cavities, as shown in Table 3. The total length of 121 2-order cavities is reduced from 3330 to 1790 2-order simplices.

**Table 3** Length distribution of 2-order cavities

| length | 8  | 10 | 12 | 14 | 16 | 18 | 20 | 22 | 24 | 26 | 28 | 30 | 32 | 34 | 36 | 38 | 1790 |
|--------|----|----|----|----|----|----|----|----|----|----|----|----|----|----|----|----|------|
| number | 27 | 28 | 15 | 12 | 3  | 8  | 6  | 6  | 2  | 1  | 1  | 5  | 1  | 1  | 3  | 2  | 121  |

The above is the result of four iterations. Each iteration is to add one 3-simplex or 2-cavity (with more than half of the same nodes) to it. Here, a 2-order cavity with length 34 is used as an example. After two iterations, the length of the cavity is reduced to 8, which achieves the best. The 2-order cavity has nodes (2,3,4,13,14,15,17,85,87,103,114,115,117,118,121,133,143,192);

another 2-order cavity with length 30 has nodes (2,3,4,13,14,15,17,85,87,103,114,115,117,121,133, 143,192), which share 27 same 2-simplices with the former, and they are deleted when two 2-order cavities are added as shown by the crossed numbers shown in Table 4. Thus the result yields a new 2-order cavity with length 10, which replace the former. The new 2-order cavity has nodes (4,13,87,117,118,133,192); another 3-simplex has nodes (4,13,87,118) to be added again. The result is a 2-order cavity with length 8, as shown in Table 4.

**Table 4** Two iterations of a 2-order cavity with length 34

| No | 0-iteration | + 2-order cavity | 1-iteration | + 3-simplex | 2-iteration |
|---|---|---|---|---|---|
| 1 | 117,118,192 | 117,133,192 | 117,118,192 | 4,13,87 | 117,118,192 |
| 2 | 4,13,118 | 4,13,87 | 4,13,118 | 4,13,118 | 13,87,133 |
| 3 | 4,87,118 | 13,15,133 | 4,87,118 | 4,87,118 | 13,117,118 |
| 4 | 13,15,133 | 13,117,133 | 13,87,133 | 13,87,118 | 87,118,192 |
| 5 | 13,87,133 |  | 13,117,118 |  | 87,133,192 |
| 6 | 13,117,118 |  | 87,118,192 |  | 117,133,192 |
| 7 | 87,118,192 |  | 87,133,192 |  | 13,117,133 |
| 8 | 87,133,192 |  | 117,133,192 |  | 13,87,118 |
| 9 |  |  | 4,13,87 |  |  |
| 10 |  |  | 13,117,133 |  |  |
| ... |  |  |  |  |  |
| 35-36 | 4,13,87 | This is a repeating simplex by formula (3) in the **Discussion** section. | | | |

The length distribution of 3-order cavities is shown in the left side of Table 5, while the shorter length is on the right side. The result is obtained by 0-1 programming[24], or through five iterations.

**Table 5** Length distribution of 3-order cavities

| length | 16 | 20 | 24 | 31 | 16 | 28 |
|---|---|---|---|---|---|---|
| number | 1 | 1 | 1 | 1 | 3 | 1 |

As can be seen above, the new method can quickly find simplices composed of all-order cavities in a sequence of nested subnetworks, but the length of the cavity is not necessarily the shortest. It needs to be shortened by multiple iterations.

**BA scale-free network**[30]

In Ref. [30], initially there are two nodes in the network. Then, one new node with two new edges are added each time according to the BA model generation mechanism. When there are 1000 nodes in the network, the process stops. Thus, a simplicial network $K = \{V, E, T\}$ is obtained, in which the Euler characteristic number (and the Betti numbers $\beta_k$) $\chi$ = 1000−1996+55 (= 1−942) = −941. The numbers of critical simplices obtained by Kannan's assignment method[30] are $c_0$=8, $c_1$=949, respectively, and the Betti persistence barcodes are obtained as shown in Fig. 5.

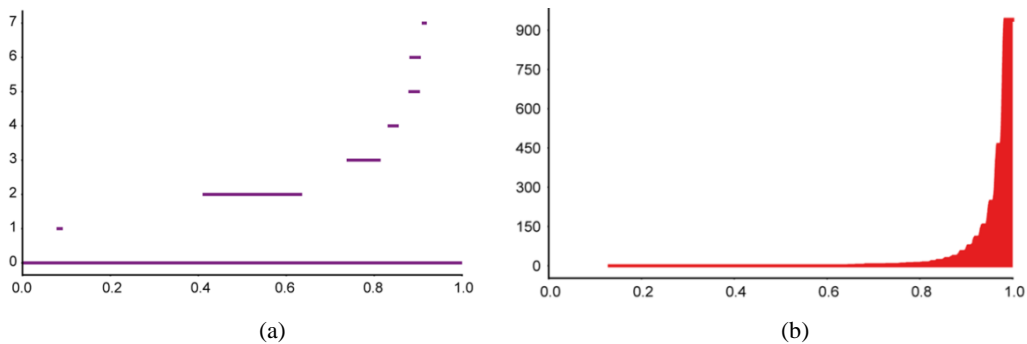

(a)    (b)

Fig. 5 (a) Betti persistence barcodes of 0-0rder cavities; (b) Betti persistence barcodes of 1-order cavities (Ref. [30])

The abscissa axis in Fig. 5 is the normalized filtration weight dividing with $w_N$, where $w_N$ is 1 plus the maximum value among weights assigned to the simplices in $K$.

In the simulation of Ref. [30], 55 2-order simplices are obtained, while our simulation result is 59 due to randomness, in which the Euler characteristic number (and the Betti numbers $\beta_k$), $\chi = 1000-1996+59$ (= $1-938$) = $-937$, and the ranks of boundary matrixes $B_1$, $B_2$ for describing the relationship between simplices are $r_1=999$, $r_2=59$, respectively. The sequences of nested simplicial subnetworks, $\varnothing \subseteq K_0 \subseteq K_1 \subseteq \cdots \subseteq K_n = K$, of this BA scale-free network $K$ are $n_1$ (= $c_0+c_1-1$ = $17+954-1$) = 970 and $n_2$ (= $r_1+\beta_k+r_2$ = $999+938+59$) = 1996, obtained by Kannan's[30] and by our assignment method, respectively. The numbers of critical simplices obtained by two assignment methods are $c_0=17$, $c_1=954$ and $c_0=1$, $c_1=938$, respectively. Obviously, our assignment method can reach the minimum value and find the representative cycles of cavities (all details are given in **Supplementary Table S3**).

**Stanford dragon graphic network**[34]
Point-cloud data of the network sample points are obtained uniformly at random from the 3-dimensional scans of the dragon photo, with reconstruction as shown in Fig. 6 (a). The point-cloud data list 1000 points in the ($x$, $y$, $z$)-coordinates. As the threshold increases from 0 to 0.016, the threshold values of all simplices in the previous ordering are recorded. The simplicial network $K$ finally formed by the point-cloud data has 1000 0-simplices (nodes), 6971 1-simplices (edges), 22712 2-simplices (triangles), and 51543 3-simplices (tetrahedrons). Its Euler characteristic number (and Betti numbers $\beta_k$) is $\chi = 1000-6971+22712-51543$ (=$1-66+1-34738$) = $-34802$. The numbers of critical simplices obtained by Otter's assignment method[34] are $c_0=1000$, $c_1=276+999=1275$, $c_2=7+210=217$, •••, and Betti persistence barcodes of 1- and 2-order cavities calculated by the method in Ref. [34] are shown in Fig. 6 (b). The lines with arrow in the picture are persistent corresponding to 1- and 2-order cavities, and the rest represent birth and death.

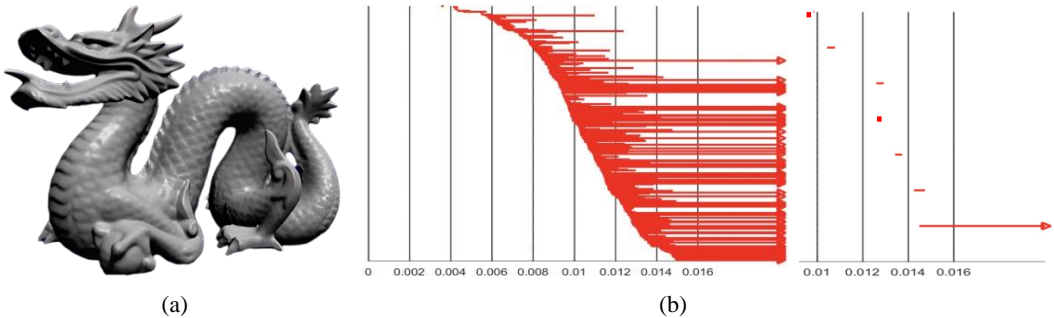

(a)        (b)

Fig. 6 (a) Reconstruction of the Stanford dragon; (b) Betti persistence barcodes of 1- and 2-order cavities (Ref. [34])

Our method can also find all the above results with noise removed, i.e. there are no barcodes of birth and death therein. The reason is that the nested sequences of the two (space and simplicial) filtrations are different. Also, the representative cycles of the only 2-order cavity have lengths 38 (Ref. [34]) and 14 (our method), respectively. All details are given in **Supplementary Table S4**.

In the above three case studies, the first example demonstrates that using spanning trees of various orders, it is possible to more easily explore higher-order network topologies from low to high; the second example shows that the network density can remain balanced as the network evolves and grows to get a filtration; the third example illustrates that increasing the threshold values will change the network from spare to dense, thereby becoming subtle and more accurate. All these examples from data science have nested simplicial subnetworks. It can be seen that, in order to study the topological characteristics of these datasets and their networks, critical simplices and Betti numbers that describe cavities are important and fundamental. Compared with the existing methods, our method can find the minimum number of critical simplices, and can also determine the number of nested simplicial subnetworks. In particular, our method can quickly find the simplices composed of all-order cavities, facilitating the understanding of the topological features of the dataset under investigation.

To further illustrate the significance of our method, consider the construction of a brain network model. Many research reports have shown that constructing the brain network model based on a fixed threshold of connectivity led to controversial results with unconvincing comparisons due to the lacks of priori biological knowledge and universally acceptable criteria. Our method does not pre-set any threshold value, but is dynamically increasing the threshold to reveal the topological characteristics of the network. This helps observe the evolution of the brain topology to obtain the best possible threshold value for the model, which clearly is related to the network order. This also helps reveal the topological characteristics of persistent homology, which is related to the number of critical simplices, thereby finding the essential differences of various brain network models. Moreover, the threshold defined in our method can be the distance between nodes, the coupling strength, or the gray level of data, etc., indicating its usefulness in a broad range of real-world applications.

## Discussion

A filtered network is a nested sequence of its subnetworks: $\varnothing \subseteq K_0 \subseteq K_1 \subseteq \cdots \subseteq K_n = K$. The filtered network is the core in topological data analysis. Therein, a filtered simplicial network induced by Morse function is particularly important. Because the persistence barcodes calculated from different nested sequences are generally different, the assigned Morse function value (or the filtration obtained by other methods) determines the number of critical simplices, the birth and death time of some cavities, and the evolution of the cavity length (from large to small).

A $k$-cavity in the homology group $Z_k/Y_k$ is a cycle in one of the linearly independent cycle-equivalent classes. Each class selects a cycle as its representative, and all representative cycles constitute a homology basis. A cavity is called optimal if the length of its representative cycle is the smallest, and a homology basis is called optimal if the total length of its representative cycles is the smallest. All the cavities of the C. elegans neural network obtained by other methods are not optimal. The optimal 1-order cavities are obtained by exhaustive searching, but the optimal 3-order cavities are obtained by the 0-1 programming method. For $k>1$, the problem of computing an optimal homology basis is NP-hard[35]. Nevertheless, the problem is polynomial-time solvable for $k=1$[36]. For an $l$-order simplicial network, without counting for the iteration of simplifying the cavity lengths, the computational complexity of the proposed method of solving equations is $O(N^3)$, where $N$ is the network size, i.e., $N = m_0+m_1+m_2+m_3+\cdots+m_l$, the total number of all-order simplices.

In the matrix equation (1), because the matrix $(T\text{-}B_k)^T(T\text{-}B_k)$ is not reversible, another method is used to solve it. To do so, consider a $k$-order orientated simplex $[i_0, i_1, \cdots, i_n]$ with a boundary $[i_0, i_1, \cdots, i_{p-1}, i_{p+1}, \cdots, i_n]$ where $i_p$ is removed. Then, each element of $B_n$ is given a plus or a minus sign[19], determined by $(-1)^p$, and the resultant matrix is denoted as $B_{[n]}$. All $k$-order cavities satisfy the following matrix equation:

$$(T\text{-}B_{[k]})(T\text{-}x)^T = (C\text{-}B_{[k]}). \tag{2}$$

Further, because the matrix $(T\text{-}B_{[k]})^T(T\text{-}B_{[k]})$ is invertible, the matrix equation (2) has a unique

solution, given by

$$(T\text{-}x)^T = [(T\text{-}B_{[k]})^T(T\text{-}B_{[k]})]^{-1}(T\text{-}B_{[k]})^T(C\text{-}B_{[k]}). \tag{3}$$

Because of different definitions of chain and cycle[19], the method would cause the problem of repeating simplices in some cycles, see the last line in Table 4. It is more appropriate to use the following formula:

$$(T\text{-}x)^T = [(T\text{-}B_{[k]})^T(T\text{-}B_{[k]})]^{-1}(T\text{-}B_{[k]})^T(C\text{-}B_{[k]}) \pmod{2}. \tag{4}$$

As warned in Ref. [34], changing the coefficient field in the definitions of chain and cycle can affect the Betti numbers. For example, if one computes the homology of the Klein bottle with coefficients in the field $F_p$, where $p$ is a primer, then $\beta_0(K) = 1$ for all primers $p$. However, $\beta_1(K) = 2$ and $\beta_2(K) = 1$ if $p=2$, but $\beta_1(K) = 1$ and $\beta_2(K) = 0$ for all other primers $p$.

With the orientated boundary matrix, the Hodge-Laplacian matrix can be obtained by $L_{(n)} = B_{[n]}^T B_{[n]} + B_{[n+1]} B_{[n+1]}^T = L_{(n)}^{down} + L_{(n)}^{up}$. It is known that the number of zero eigenvalues of the Hodge-Laplacian matrix $L_{(k)}$ equals the number of $k$-order cavities, namely the Betti number $\beta_k$. And note that zero eigenvectors can provide information about cavities, which is a topic worth studying.

**Acknowledgements**

The authors would like to thank the research support of the National Natural Science Foundation of China under Grants No. 62173095 and No. 12005001, and by the Hong Kong Research Grants Council through General Research Funds under Grant CityU 11206320.

**Author contributions**

DS and GC wrote the text and developed the theory. CM proposed another method to solve the matrix equation. ZC and CM performed the simulations and computations for cross check. All authors checked and finalized the entire manuscript.

**Competing interests**

Authors declare no competing interests.

**Supplementary information** is provided in the Appendix.

**Appendix: Supplementary Information (Tables S1 to S4)**

**Table S1: Torustriangulation network**

This table lists the Morse function values of all simplices, the spanning trees of the boundary matrixes $B_1$ and $B_2$, and the simplices composed of 1- and 2-order cavities for the network.

**Table S2: C. elegans neural network[33]**

This table lists the Morse function values of all simplices, the simplices composed of 1-, 2- and 3-order cavities, the iterative process of 2-order cavities, and an iterative example for the network.

**Table S3: BA scale-free model network[30]**

This table lists the Morse function values of all simplices, the simplices composed of 1-order cavities, and the results obtained by Kannan's[30] method for the network which is simulated here.

Variables in Kannan's method are as follows:

DFM=Discrete Morse function values;

Flag=To keep track with the size of the set $U_\alpha$ for each simplex $\alpha$;

IsCritical=To indicate if a given simplex is critical;

FiltrationWeight=To store the filtration weight corresponding to each simplex.

**Table S4: Stanford dragon graphicnetwork[34]**

This table lists the points in the ($x$, $y$, $z$)-coordinates, the present thresholds of all simplices, the persistence barcodes of 1- and 2-order cavities calculated by **javaplex** in Ref. [34], the representative cycles with two lengths of the only 2-order cavity, and the Morse function values of all simplices obtained by the new method for the network.

Data of **Tables S1 to S4** are available:

https://github.com/ChuangMa1900/Supplementary-Information-Tables-S1-to-S4.git